\documentclass[11pt,a4paper]{amsart}
\usepackage[margin=27mm,headheight=15pt]{geometry}
\usepackage[T1]{fontenc}
\usepackage{lmodern,amsmath,amssymb,amsthm,mathrsfs,microtype}
\usepackage{enumitem,tikz}
\usetikzlibrary{arrows.meta}
\usepackage[colorlinks=true,linkcolor=blue,citecolor=blue,urlcolor=blue]{hyperref}
\usepackage{bookmark}
\hypersetup{pdftitle={Prime-level congruence rigidity for SL2(Fq[t])},
  pdfauthor={Sanghoon Kwon}}
\setlist[enumerate]{label=(\roman*),itemsep=2pt,topsep=4pt}
\newtheorem{theorem}{Theorem}[section]
\newtheorem{lemma}[theorem]{Lemma}
\newtheorem{proposition}[theorem]{Proposition}
\newtheorem{corollary}[theorem]{Corollary}
\newtheorem{mainthm}{Theorem}

\newtheorem{maincor}{Corollary}

\theoremstyle{definition}
\newtheorem{example}[theorem]{Example}
\newcommand{\F}{\mathbb F}
\newcommand{\SL}{\operatorname{SL}_2}
\newcommand{\GL}{\operatorname{GL}_2}
\newcommand{\PSL}{\operatorname{PSL}_2}
\newcommand{\Aut}{\operatorname{Aut}}
\newcommand{\Inn}{\operatorname{Inn}}
\newcommand{\Good}{\operatorname{Good}}
\newcommand{\Target}{\operatorname{Target}}
\newcommand{\ev}{\operatorname{ev}}
\newcommand{\diag}{\operatorname{diag}}
\newcommand{\ql}{\operatorname{ql}}
\newcommand{\cl}{\operatorname{cl}_{\rm cong}}
\newcommand{\DF}{\mathfrak D}
\title[Prime-level congruence rigidity]{Prime-level congruence rigidity for
  $\mathrm{SL}_2(\mathbb F_q[t])$}
\author[S. Kwon]{Sanghoon Kwon}
\address{Department of Mathematics Education, Catholic Kwandong University,
24 Beomil-ro 579beon-gil, Gangneung 25601, Republic of Korea}
\email{skwon@cku.ac.kr}
\subjclass[2020]{Primary 20E36; Secondary 11F06, 20E08, 20G30}
\keywords{Drinfeld modular group, congruence subgroup, Reiner automorphism, tree lattice}
\date{}
\begin{document}

\begin{abstract}
Let $A=\mathbb F_q[t]$. We classify the automorphisms of
$\mathrm{SL}_2(A)$ that send the principal congruence subgroup $\Gamma(p)$
to a congruence subgroup for infinitely many primes $p$. They are precisely
those obtained by applying a ring automorphism of $A$ entrywise and then
conjugating by a matrix in $\mathrm{GL}_2(A)$. For any finite set of
primes, we construct a Reiner automorphism sending $\Gamma(p)$ to a
congruence subgroup at exactly those primes. Let $L$ be the linear map
defining a Reiner automorphism. For each prime $s$, $\Gamma(s)$ is the
image of some prime-level subgroup $\Gamma(p)$ precisely when
$f\mapsto L(f)\bmod s$ is multiplicative on $A$. If the residue field
at $p$ has at least four elements, then every noncongruence
image of $\Gamma(p)$ has the same reduction as $\mathrm{SL}_2(A)$ modulo
every nonzero polynomial.
\end{abstract}
\maketitle
\tableofcontents

\section{Introduction}
Let $k=\F_q$, $A=k[t]$, and $\Gamma=\SL(A)$. For $0\ne a\in A$, the
matrices in $\Gamma$ congruent to the identity modulo $a$ form the
principal congruence subgroup $\Gamma(a)$. A subgroup is called
congruence if it contains some $\Gamma(a)$. Throughout, a prime
means a monic irreducible polynomial in $A$.

Fix an automorphism $\alpha$ of $\Gamma$ and let the prime $p$ vary.
The quotient of $\Gamma$ by $\alpha(\Gamma(p))$ is always isomorphic
to $\SL(A/(p))$. The subgroup $\alpha(\Gamma(p))$ itself can be
noncongruence. If it is congruence for infinitely many $p$, then
$\alpha$ is a composition of an entrywise ring automorphism of $A$
and conjugation in $\GL(A)$ (Corollary~\ref{cor:main}).

For the automorphisms introduced by Reiner, Theorem~\ref{thm:main}
determines exactly when the image of $\Gamma(p)$ is a congruence subgroup.
Put $S=\SL(k)$ and
$u(f)=\left(\begin{smallmatrix}1&f\\0&1\end{smallmatrix}\right)$
for $f\in A$. Let $\mathscr R_q$ be the set of $k$-linear bijections
of $A$ fixing $k$ pointwise. For $L\in\mathscr R_q$, the associated
\emph{Reiner automorphism} $\Phi_L$ is the unique automorphism of
$\Gamma$ fixing $S$ pointwise and sending $u(f)$ to $u(L(f))$.
Existence and uniqueness follow from Nagao's decomposition~\eqref{eq:amalgam}.

\begin{samepage}
To test multiplicativity modulo a prime, we use the
\emph{multiplicative defects} of $L$ and the ideal they generate.
\begin{equation}\label{eq:defect}
 \Delta_L(f,g)=L(fg)-L(f)L(g),\qquad
 \DF_L=(\Delta_L(f,g)\mid f,g\in A)=(d_L).
\end{equation}
We take $d_L$ monic for a nonzero ideal and put $d_L=0$ otherwise.
Every prime divides zero.
\end{samepage}

Call a prime $p$ a \emph{good source prime} if $\Phi_L(\Gamma(p))$ is
congruence. A \emph{target prime} is a prime $s$ for which
$\Phi_L(\Gamma(p))=\Gamma(s)$ for some prime $p$. Write $\Good(L)$ and
$\Target(L)$ for the respective sets. We call $L$ \emph{degree-preserving}
if it preserves each subspace $A_{\le n}$ of polynomials of degree at
most $n$, including zero. This extra condition enters the constructions
in Sections~\ref{sec:finite} and~\ref{sec:tree}.

\begin{mainthm}\label{thm:main}
Let $L\in\mathscr R_q$. Every good source prime $p$ has a unique target
$s_p$ satisfying $\Phi_L(\Gamma(p))=\Gamma(s_p)$. The correspondence is
a bijection
\[
 \Good(L)\longrightarrow\Target(L)=\{s\text{ prime}:s\mid d_L\},
 \qquad \deg p=\deg s_p .
\]
For a target $s$, its source $p$ is the minimal polynomial over $k$ of
$L(t)\bmod s$. If $L$ is not multiplicative, then for every nonzero
defect $D=\Delta_L(f,g)$,
\[
 \sum_{p\in\Good(L)}\deg p=\deg\operatorname{rad}(d_L)\le\deg D .
\]
Here $\operatorname{rad}(d_L)$ is the product of the distinct monic prime
divisors, with empty product $1$. Every prescribed finite set occurs as
$\Good(L)$ for a degree-preserving $L$ attaining equality for some $D$.

The following conditions are equivalent.
\begin{enumerate}
\item $\Good(L)$ is infinite.
\item $\DF_L=0$.
\item $L(f)(t)=f(at+b)$ for some $a\in k^\times$ and $b\in k$.
\item $\Phi_L$ sends every congruence subgroup to a congruence subgroup.
\end{enumerate}
\end{mainthm}

Call $\alpha$ \emph{standard} if $\alpha=\operatorname{Ad}_g\circ\sigma_*$,
where $g\in\GL(A)$ acts by conjugation and $\sigma_*$ applies a ring
automorphism $\sigma$ of $A$ entrywise. Such a $\sigma$ acts on $k$ by a
field automorphism and sends $t$ to $at+b$, $a\in k^\times$, $b\in k$.
For example, the substitution $t\mapsto t+1$ gives a standard Reiner
automorphism sending $\Gamma(t)$ to $\Gamma(t+1)$.

\begin{maincor}\label{cor:main}
For every prime power $q$ and every $\alpha\in\Aut(\Gamma)$, the following
are equivalent.
\begin{enumerate}
\item $\alpha(\Gamma(p))$ is congruence for infinitely many primes $p$.
\item $\alpha$ sends every congruence subgroup to a congruence subgroup.
\item $\alpha$ is standard.
\end{enumerate}
For $\alpha=\operatorname{Ad}_g\circ\sigma_*$ we have
$\alpha(\Gamma(a))=\Gamma(\sigma(a))$ for every $a\ne0$, and
$\alpha^{-1}$ also preserves congruence subgroups.
\end{maincor}

Write $\rho_a:\Gamma\to\SL(A/(a))$ for reduction modulo $a\ne0$.
A finite quotient of $\Gamma$ is a congruence quotient if its quotient
map factors through some $\rho_a$. The cosets of the $\Gamma(a)$ form a basis for the
congruence topology. For any finite list of congruence quotients, there
is a nonstandard automorphism inducing the identity on them all
(Corollary~\ref{cor:tests}).

For a normal subgroup $N$ of $\Gamma$, its \emph{quasi-level}
$\ql(N)=\{f\in A:u(f)\in N\}$ records its upper unipotent elements.
This is an additive subgroup of $A$. The largest ideal of $A$ contained
in it is the \emph{level} of $N$ \cite[\S4]{MS12}. Thus $\ql(\Gamma(p))=pA$
and $\ql(\Phi_L(\Gamma(p)))=L(pA)$. These quasi-levels can agree even
when $\Phi_L(\Gamma(p))$ is noncongruence.

\begin{example}\label{ex:same}
Take $k=\F_3$ and $p=t^2+1$, which has no root in $k$.
Write $f=a+bt+pg$ with $a,b\in k$ and $g\in A$, and set $L(f)=f+b$.
Subtracting $b$ gives the inverse. Thus $L$ preserves degrees and fixes
$pA$ pointwise, giving $\ql(\Phi_L(\Gamma(p)))=pA$. But $L(t)=t+1$,
$L(t^2)=t^2$ and $L(pt)=pt$, whence
\[
   \Delta_L(t,t)=t-1,\qquad \Delta_L(p,t)=-p.
\]
The two defects are relatively prime. Hence $d_L=1$ and $\Good(L)$ is empty.
In fact, the image of $\Gamma(p)$ surjects onto every finite congruence
quotient, since $|A/(p)|=9$
(Corollary~\ref{cor:density}).
\end{example}

Reiner introduced this construction for matrices over a polynomial ring
\cite{Reiner}. Mason and Schweizer extended it to rings of functions
regular away from one fixed place and constructed noncongruence subgroups
of prescribed level \cite{MS06,MS12}. For the terminology see
\cite[Remark~3.7 and Definition~3.9]{MS26}. In their proof of
\cite[Theorem~3.10]{MS26} the level is fixed and the automorphism varies.
Here we fix the automorphism and determine all its good source primes.
Mason and Schweizer call an automorphism standard if it preserves
congruence subgroups \cite{MS16}. Corollary~\ref{cor:main} identifies
this condition with our algebraic definition for $k[t]$.

Section~\ref{sec:tree} determines the congruence closures of the
automorphic images of $\Gamma(t)$. In the construction of Mason and
Schweizer \cite{MS06}, these subgroups are indexed by $k$-linear forms
$\ell:A\to k$ with $\ell(1)=1$. The associated hyperplane is
$\ker\ell\subset A$. All these subgroups have finite index in $\Gamma$
and hence are arithmetic. Their tree actions are conjugate, and the
vertex-averaging operators on their quotient graphs are unitarily
equivalent. Their congruence closures can differ. For each $q\in\{2,3\}$, exactly
six noncongruence members have proper closure. Each of these closures
has index $q$ in $\Gamma$ (Proposition~\ref{prop:small}). For background on tree
lattices and the congruence subgroup problem, see \cite{BL,Lub91,MPSZ}.

\section{Kernels of reduction modulo a prime}\label{sec:recognition}
Write $B(A)$ for the upper triangular subgroup of $\Gamma$ and
$B(k)=B(A)\cap S$. For a unit $a$ write $h(a)=\diag(a,a^{-1})$.
In Nagao's decomposition of $\GL(A)$ \cite{Nagao,Serre}, the two factors
share the matrices $\diag(c,1)$ for $c\in k^\times$. Taking the kernel
of the determinant thus gives the amalgamated product
\cite[\S5.1 of the preprint]{MitraSankaran}
\begin{equation}\label{eq:amalgam}
 \Gamma=S*_{B(k)}B(A).
\end{equation}
For $L\in\mathscr R_q$, applying $L$ to the upper-right entry defines
an automorphism of $B(A)$ by $k$-linearity. It fixes $B(k)$ pointwise
and hence extends by \eqref{eq:amalgam} to $\Phi_L$, with inverse
$\Phi_{L^{-1}}$. We also use that $\Gamma$ is generated by upper and lower
elementary unipotent matrices, by the Euclidean algorithm in $k[t]$.
The lower ones are conjugate to upper ones by
$w=\left(\begin{smallmatrix}0&1\\-1&0\end{smallmatrix}\right)\in S$.

For a prime $s$ put $\kappa(s)=A/(s)$ and denote the projective reduction
$\Gamma\to\PSL(\kappa(s))$ by $\bar\rho_s$. We use the same matrices
$u(x)$, $w$ and $h(a)$ over any finite field, and write
$u(E)=\{u(x):x\in E\}$ for sets of parameters.

\begin{lemma}\label{lem:kernel}
Let $R=q^d$, $d\ge1$, and put $T_R=\SL(\F_R)$ and $P_R=T_R/Z(T_R)$.
\begin{enumerate}
\item If a surjection $\Gamma\to P_R$ factors through a finite congruence
quotient, its kernel is $\ker\bar\rho_s$ for a unique prime $s$ of degree $d$.
\item Suppose $\phi:\Gamma\twoheadrightarrow T_R$ and its composition
$\bar\phi:\Gamma\to P_R$ factors through a finite congruence quotient.
Then $\ker\phi=\Gamma(s)$ for the same prime $s$. In particular, every
normal congruence subgroup with quotient isomorphic to $T_R$ has this form.
\end{enumerate}
\end{lemma}
\begin{proof}
(i) Choose $a=\prod_i s_i^{e_i}$ with $\Gamma(a)$ in the kernel, the
$s_i$ being distinct primes. Reduction gives a surjection
\[
 \Gamma \longrightarrow \SL(A/(a))
      \simeq \prod_i\SL(A/(s_i^{e_i}))\longrightarrow P_R .
\]
Indeed elementary matrices generate over each finite local ring by row
elimination, and their entries lift one component at a time by the Chinese
remainder theorem. Let $p_0=\operatorname{char}k$. The kernel of reduction
from the product to $\prod_i\SL(\kappa(s_i))$ is a $p_0$-group, since its
successive congruence quotients embed in additive matrix groups.
Its image in $P_R$ is trivial. For $R\ge4$ this follows from the
nonabelian simplicity of $P_R$. In the two remaining cases,
$P_2\simeq S_3$ has no nontrivial normal $2$-subgroup, and
$P_3\simeq A_4$ has no nontrivial normal $3$-subgroup.

The factors $\SL(\kappa(s_i))$ have commuting normal images. For $R\ge4$,
simplicity forces one factor to surject and the others to have trivial image. Write
$Q=|\kappa(s_i)|$ for the size of this factor's residue field. The
groups $T_2$ and $T_3$ are solvable, so $Q\ge4$. The group $T_Q$ is
perfect and its proper normal subgroups are central. Indeed a noncentral
normal subgroup surjects onto the simple group $P_Q$ and hence has
abelian quotient in $T_Q$. The kernel of $T_Q\twoheadrightarrow P_R$ is
therefore central and has order prime to $p_0$. The largest powers of $p_0$
dividing $|T_Q|$ and $|P_R|$ are $Q$ and $R$, respectively, so $Q=R$.
Comparison of the orders gives precisely $Z(T_Q)$ as kernel.

For $R=2,3$ the target is solvable, so the perfect factors with $Q\ge4$
have trivial image. Here $q=R$. The remaining factors are copies of
$\SL(k)$ indexed by degree-one primes $t-c$. The proper normal subgroups
of $S_3$ lie in $A_3$, and those of $A_4$ lie in the Klein four group.
Thus one factor must surject. The other images lie in the trivial center
of $P_R$. For $R=2$ the resulting surjection is an isomorphism.
For $R=3$ its kernel is the center of order two, since $-I$ is the
only involution in $\SL(\F_3)$. This proves
the asserted form of the kernel in all cases. Uniqueness follows from
$\ker\bar\rho_s\cap u(A)=u(sA)$.

(ii) Identify $\kappa(s)$ with $\F_R$. By (i), $\bar\phi$ and
$\bar\rho_s$ differ by an automorphism of $P_R$. For $R\ge4$ all such
automorphisms come from conjugation in $\GL(\F_R)$ and field automorphisms
\cite[\S2]{White}, and hence lift to $T_R$. For $R=2$ use $P_2=T_2$.
For $R=3$, $\Aut(A_4)\simeq S_4$ is induced by
$\mathrm{PGL}_2(3)\simeq S_4$ acting by conjugation.
Choose a lift $\beta\in\Aut(T_R)$. The map
$\gamma\mapsto\phi(\gamma)\bigl(\beta\rho_s(\gamma)\bigr)^{-1}$
takes values in $Z(T_R)$ and is a homomorphism by centrality. In even
characteristic this center is trivial. In odd characteristic it has
order two, and the elementary unipotents generating $\Gamma$ have odd
order whenever they are nonidentity. Thus $\phi=\beta\circ\rho_s$.
\end{proof}

The \emph{congruence closure} of $N\le\Gamma$ is
$\cl(N)=\bigcap_{a\ne0}N\Gamma(a)$, where the intersection is over
nonzero polynomials. It equals $\Gamma$ precisely when $N$ surjects onto
every finite congruence quotient. We then call $N$ \emph{congruence-dense}.

\begin{corollary}[Congruence closure at prime levels]\label{cor:density}
Suppose $N\lhd\Gamma$ and $\Gamma/N\simeq\SL(\F_{q^d})$ with $q^d\ge4$.
Then either $N=\Gamma(s)$ for a unique prime $s$ of degree $d$, or
$N$ is congruence-dense. In particular, for any $\alpha\in\Aut(\Gamma)$
and any prime $p$ with $q^{\deg p}\ge4$, a noncongruence image
$\alpha(\Gamma(p))$ is congruence-dense.
\end{corollary}
\begin{proof}
Take $\phi:\Gamma\twoheadrightarrow T_{q^d}$ with kernel $N$. If the
closure is proper, then $N\Gamma(a)\ne\Gamma$ for some $a$. Its image
$\phi(\Gamma(a))$ is a proper normal subgroup of $T_{q^d}$, hence central.
Thus $\bar\phi$ factors through reduction modulo $a$ and
Lemma~\ref{lem:kernel}(ii) applies. Conversely each $\Gamma(s)$ is closed
and proper. For the last assertion, take $N=\alpha(\Gamma(p))$.
\end{proof}

\begin{lemma}\label{lem:field}
Let $F,F'$ be finite fields. Suppose an isomorphism
$\beta:\SL(F)\to\SL(F')$ sends $w$ to $w$ and satisfies
$\beta(u(x))=u(\theta(x))$ for an additive bijection $\theta:F\to F'$
with $\theta(1)=1$. Then $\theta$ is a field isomorphism.
\end{lemma}
\begin{proof}
The normalizers of $u(F)$ and $wu(F)w^{-1}$ are the stabilizers of their
unique fixed points $\infty$ and $0$ on $\mathbb P^1(F)$. Their intersection
is the diagonal group, so $\beta$ sends diagonals to diagonals.
For $a\ne0$ put $x=\theta(a)$ and $y=\theta(a^{-1})$.
Applying $\beta$ to $h(a)=u(a)w\,u(a^{-1})w\,u(a)w$ gives
\[
 \beta(h(a))=\begin{pmatrix}2x-x^2y&xy-1\\1-xy&y\end{pmatrix}.
\]
Diagonality gives $xy=1$, hence $\beta(h(a))=h(\theta(a))$.
Since $h(ab)=h(a)h(b)$ and $h$ is injective, $\theta(ab)=\theta(a)\theta(b)$
for nonzero $a,b$. Additivity gives $\theta(0)=0$.
\end{proof}

\begin{figure}[htbp]
\centering
\begin{tikzpicture}[>=Stealth,font=\small]
\node (G) at (0,0) {$\Gamma$};
\node (Gp) at (5.8,0) {$\Gamma$};
\node (F) at (0,-1.15) {$\SL(\kappa(p))$};
\node (Fp) at (5.8,-1.15) {$\SL(\kappa(s))$};
\draw[->] (G)--node[above]{$\Phi_L$}(Gp);
\draw[->] (G)--node[left]{$\rho_p$}(F);
\draw[->] (Gp)--node[right]{$\rho_s$}(Fp);
\draw[->] (F)--node[above]{$\beta_p$}(Fp);
\end{tikzpicture}
\caption{For $\Phi_L(\Gamma(p))=\Gamma(s)$, the group isomorphism $\beta_p$
sends $u(f\bmod p)$ to $u(L(f)\bmod s)$. The resulting map on entries is
the field isomorphism $\theta_p:\kappa(p)\to\kappa(s)$.}
\label{fig:residue}
\end{figure}
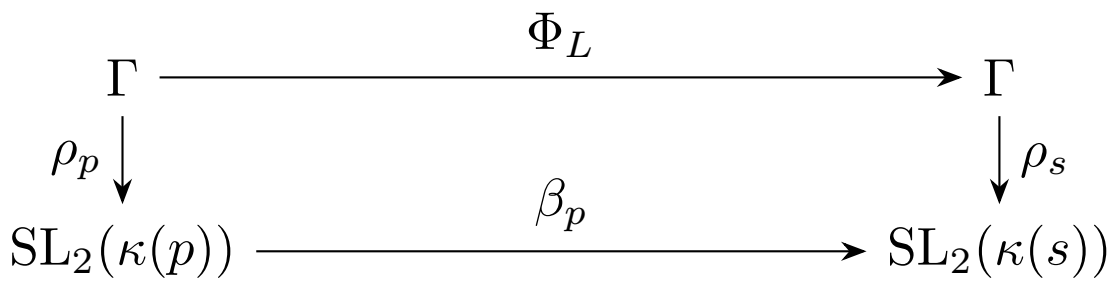

\begin{proposition}\label{prop:residue}
For $p\in\Good(L)$ there is a unique prime $s$ of the same degree with
$\Phi_L(\Gamma(p))=\Gamma(s)$. The induced map
$\theta_p:\kappa(p)\to\kappa(s)$, $f\bmod p\mapsto L(f)\bmod s$,
is a $k$-algebra isomorphism.
\end{proposition}
\begin{proof}
The image is a normal congruence subgroup with quotient isomorphic to
$\SL(\kappa(p))$, so Lemma~\ref{lem:kernel} gives $s$.
Intersecting with $u(A)$ gives $L(pA)=sA$, so $\theta_p$ is a well-defined
additive bijection. Let $\beta_p:\SL(\kappa(p))\to\SL(\kappa(s))$ be the
induced group isomorphism, characterized by
$\beta_p\circ\rho_p=\rho_s\circ\Phi_L$ (Figure~\ref{fig:residue}).
It fixes $w$ and sends $u(x)$ to $u(\theta_p(x))$ for $x\in\kappa(p)$.
Lemma~\ref{lem:field} gives multiplicativity, and $L$ fixes $k$.
\end{proof}

\section{The defect ideal and rigidity}\label{sec:rigidity}
The defects $\Delta_L(t,t^n)$ for $n\ge1$ generate the ideal in \eqref{eq:defect}.
For $h=L(t)$ we have
\begin{equation}\label{eq:generator}
 \DF_L=(L(t^n)-h^n\mid n\ge2)
      =(\Delta_L(t,t^n)\mid n\ge1).
\end{equation}
Indeed, modulo the first ideal, $k$-linearity gives $L(f)\equiv f(h)$,
so all multiplicative defects vanish. Conversely $E_n=L(t^n)-h^n$
satisfies $E_1=0$ and $E_{n+1}=\Delta_L(t,t^n)+hE_n$. Thus each $E_n$
belongs to the last ideal, which lies in $\DF_L$ by definition.
In particular, $b\mid d_L$ precisely when $f\mapsto L(f)\bmod b$ is
multiplicative, for every nonconstant $b\in A$.

\begin{proof}[Proof of Theorem~\ref{thm:main}]
For a good source $p$, Proposition~\ref{prop:residue} gives
$L(f)\bmod s_p=\theta_p(f\bmod p)$ for all $f\in A$. Since
$\theta_p:\kappa(p)\to\kappa(s_p)$ is a field isomorphism, $s_p$ divides
every $\Delta_L(f,g)$ and hence $d_L$.
Distinct sources have distinct targets by invertibility of $\Phi_L$ and
the identity $\Gamma(p)\cap u(A)=u(pA)$.

Conversely, suppose $s\mid d_L$. Then $\psi_s(f)=L(f)\bmod s$ is a
surjective unital $k$-algebra map $A\to\kappa(s)$, by the surjectivity
of $L$. Hence $\psi_s(f)=f(\bar h)$ for $\bar h=L(t)\bmod s$ and
$k[\bar h]=\kappa(s)$. Its kernel is generated by the monic minimal
polynomial $p$ of $\bar h$, giving $\deg p=\deg s$. Applying $\psi_s$
entrywise gives the same homomorphism as $\rho_s\Phi_L$ on $S$ and
$B(A)$, hence on $\Gamma$ by \eqref{eq:amalgam}. Its kernel is
$\Gamma(p)$, so $\Phi_L(\Gamma(p))=\Gamma(s)$. This proves the bijection.
Applying it to the inverse gives the corresponding source formula
\begin{equation}\label{eq:inverse}
 \Good(L)=\Target(L^{-1})=\{p\text{ prime}:p\mid d_{L^{-1}}\}.
\end{equation}

If $d_L\ne0$, the number of good sources is the number of distinct prime
divisors of $d_L$. Their degree sum is $\deg\operatorname{rad}(d_L)$,
since the bijection preserves residue degrees. The radical divides every
nonzero defect $D$, giving the stated bound. The construction in
Proposition~\ref{prop:exact} attains it for every finite set.

Infinitely many good primes therefore force $d_L=0$. In that case $L$
is a $k$-algebra automorphism. Writing $h=L(t)$, surjectivity gives
$F(h)=t$ for some polynomial $F$, so $\deg F\deg h=1$ and $h=at+b$
with $a\ne0$. Conversely, for any such substitution $\Phi_L$ acts
entrywise and sends $\Gamma(a_0)$ to $\Gamma(L(a_0))$ for every
$a_0\ne0$. It therefore preserves congruence subgroups, which in turn
implies that every prime is good.
\end{proof}

We reduce an arbitrary automorphism to one fixing $S$ pointwise and
preserving $u(A)$ by composing it with a standard automorphism,
as in \cite{MitraSankaran}. For $q=2$ see also
\cite[Proposition~5.2]{MS26}.

\begin{lemma}\label{lem:normalization}
Every $\alpha\in\Aut(\Gamma)$ can be written $\delta\circ\Phi_L$, where
$\delta$ is standard and $L\in\mathscr R_q$.
\end{lemma}
\begin{proof}
Let $T$ be the Bruhat--Tits tree for $\SL$ over $k(\!(t^{-1})\!)$,
the completion of $k(t)$ at infinity. The quotient $\Gamma\backslash T$
is a ray, with depth measured from its initial vertex. Its standard lift
to $T$ meets each $\Gamma$-orbit of vertices or edges exactly once
(Section~\ref{sec:tree}). Since $\Gamma$ preserves the two vertex classes
of the bipartite tree, the finite subgroup $\alpha(S)$ fixes a vertex.
Along this lift the stabilizer is $S$ at depth zero. At positive depth
it is upper triangular with a normal $p_0$-subgroup and quotient
$k^\times$, where $p_0=\operatorname{char}k$ \cite{Serre}. The norm-one
subgroup of $\F_{q^2}^\times$ embeds in $S$ by multiplication on the
two-dimensional space $\F_{q^2}$ over $k$. Its order $q+1$ is prime
to $p_0$ and exceeds $q-1$. Thus $\alpha(S)$ is contained in a
conjugate of $S$ and equals it by order. After an inner conjugation
assume $\alpha(S)=S$.

The automorphism induced on $S/Z(S)$ comes from a matrix conjugation
and a field automorphism, as in Lemma~\ref{lem:kernel}(ii)
(for $q=2$, all automorphisms of $S_3$ are inner). These operations extend
to $\Gamma$, extending the coefficient automorphism by fixing $t$.
Their action on $S$ agrees with $\alpha|_S$, since any discrepancy is a
homomorphism $S\to Z(S)$ and elementary unipotents rule it out as above.
Composing with the inverse standard map makes $\alpha$ fix $S$ pointwise.

Direct calculation gives $C_\Gamma(u(1))=Z(\Gamma)u(A)$, where
$Z(\Gamma)=\{\lambda I:\lambda\in k^\times,\ \lambda^2=1\}$.
The elements of $p_0$-power order in this centralizer form exactly $u(A)$.
Thus $\alpha(u(f))=u(L(f))$ for an additive bijection $L$ fixing $k$.
Conjugation by $h(a)$ gives $L(a^2 f)=a^2 L(f)$ for every $a\in k^\times$.
In even characteristic each scalar is a square. In odd characteristic
each $c\in k$ is the difference of the squares of $(c+1)/2$ and $(c-1)/2$.
Hence $L$ is $k$-linear. Since $S$ and $u(A)$ generate $\Gamma$,
the normalized automorphism is $\Phi_L$.
\end{proof}

\begin{proof}[Proof of Corollary~\ref{cor:main}]
Conjugation by $g\in\GL(A)$ fixes $\Gamma(a)$, and $\sigma_*$ sends it
to $\Gamma(\sigma(a))$. Thus standard maps preserve congruence in both
directions. Write $\alpha=\delta\circ\Phi_L$ as in
Lemma~\ref{lem:normalization}. Infinitely many congruence images for
$\alpha$ give as many for $\Phi_L$, so Theorem~\ref{thm:main} makes $L$
an affine substitution and $\alpha$ standard. Finally, preservation of
all congruence subgroups implies condition~(i) of Corollary~\ref{cor:main}.
\end{proof}

For $q=2$ the group of standard automorphisms of $\Gamma$ is
$\Inn(\Gamma)\rtimes\langle t\mapsto t+1\rangle$. Here $\GL(A)=\Gamma$,
and translation interchanges the distinct normal subgroups $\Gamma(t)$
and $\Gamma(t+1)$. For odd $q$, conjugation by $\diag(\nu,1)$ with
$\nu\in k^\times$ a nonsquare is outer. Indeed any matrix centralizing
both $u(k)$ and $w u(k)w^{-1}$ is scalar, and a scalar matrix
has square determinant.

\section{Prescribed finite sets of good primes}\label{sec:finite}
Suppose $\deg L(t)=1$, as holds in particular when $L$ preserves degrees.
Let $\sigma$ be the $k$-algebra automorphism with $\sigma(t)=L(t)=at+b$.
The normalized map $M=\sigma^{-1}\circ L$ fixes $t$. Write
$J_M=(M(f)-f\mid f\in A)=(j_M)$ with the same monic convention as for $d_L$.

\begin{proposition}\label{prop:exact}
We have $\DF_L=\sigma(J_M)$. A prime $p$ is good precisely when
$p\mid j_M$, and then $\Phi_L(\Gamma(p))=\Gamma(\sigma(p))$ with monic
target $s_p=a^{-\deg p}\sigma(p)$. For every finite set $\mathcal P$ of
primes there is a degree-preserving $M$ fixing $t$ for which
$\Good(M)=\mathcal P$ and equality holds in the degree bound of
Theorem~\ref{thm:main}.
\end{proposition}
\begin{proof}
By \eqref{eq:generator} and $k$-linearity, $\DF_M=J_M$.
Applying $\sigma$ to each defect gives $\DF_L=\sigma(J_M)$.
If $p\mid j_M$, then $M(f)\equiv f\bmod p$ for every $f$, so
$\rho_p\Phi_M=\rho_p$ on both factors of \eqref{eq:amalgam} and
$\Phi_M(\Gamma(p))=\Gamma(p)$. Conversely, for $p\in\Good(M)$,
Proposition~\ref{prop:residue} gives a $k$-isomorphism
$\kappa(p)\to\kappa(s)$ taking $t\bmod p$ to $t\bmod s$, since $M(t)=t$.
Their monic minimal polynomials are $p$ and $s$, so $p=s$. Since
$t\bmod p$ generates $\kappa(p)$ over $k$, the isomorphism is the
identity, giving $p\mid j_M$. As $\sigma_*$ permutes the principal
congruence subgroups, $\Phi_L=\sigma_*\circ\Phi_M$ has the same good
sources as $\Phi_M$ and sends $\Gamma(p)$ to $\Gamma(\sigma(p))$ at
each such prime.

Take $r=\prod_{p\in\mathcal P}p$, using $r=1$ for the empty set, and
choose $m>\max\{1,\deg r\}$. Define on the monomial basis
\[
 M(t^m)=t^m+r,\qquad M(t^n)=t^n\quad(n\ne m).
\]
Its inverse subtracts $r$ from $t^m$, so $M$ preserves degrees and
$J_M=(r)$. Moreover $\Delta_M(t,t^{m-1})=r\ne0$ and
$\deg r=\sum_{p\in\mathcal P}\deg p$, as required.
\end{proof}

\begin{example}\label{ex:two}
Over $\F_2$, the choice $r=t(t+1)$ and $m=3$ gives
$\Good(M)=\{t,t+1\}$ and $\Delta_M(t,t^2)=t(t+1)$.
To see distinct source and target primes, take instead $r=t$, $m=2$
and $L=\sigma\circ M$ with $\sigma(t)=t+1$. Then $\Good(L)=\{t\}$,
$\Target(L)=\{t+1\}$ and $\Delta_L(t,t)=t+1$. In particular,
$\Phi_L(\Gamma(t))=\Gamma(t+1)$.
\end{example}

\begin{corollary}\label{cor:tests}
For every $a\ne0$ there is a degree-preserving $M\in\mathscr R_q$ with
finite $\Good(M)$ such that $\Phi_M$ induces the identity on
$\Gamma/\Gamma(a)$.
\end{corollary}
\begin{proof}
Take $r=a$ and $m>\max\{1,\deg a\}$ in the construction above. Then
$(M-\mathrm{id})(A)\subseteq aA$, so $\rho_a\Phi_M=\rho_a$ on both
amalgam factors. The nonzero defect $\Delta_M(t,t^{m-1})=a$ gives
finiteness. For several congruence quotients take the product of
moduli through which their maps factor.
\end{proof}

\section{Congruence closures in the orbit of \texorpdfstring{$\Gamma(t)$}{Gamma(t)}}\label{sec:tree}
Let $\ell:A\to k$ be $k$-linear with $\ell(1)=1$, so $A=k\oplus\ker\ell$.
Following Mason and Schweizer \cite[Theorem~3.5 and Lemma~4.4]{MS06}, fix
$S$ pointwise and replace the upper-right entry $f$ in $B(A)$ by $\ell(f)$.
These rules agree on $B(k)$ and define a retraction $r_\ell:\Gamma\to S$
by \eqref{eq:amalgam}. Its kernel $H_\ell$ has quasi-level $\ker\ell$.
For evaluation $\ev_c(f)=f(c)$ at $c\in k$, we have
$\ker\ev_c=(t-c)A$ and $H_{\ev_c}=\Gamma(t-c)$.

The vertices of the Bruhat--Tits tree $T$ at infinity are homothety
classes of $\mathcal O_\infty$-lattices in $k(\!(t^{-1})\!)^2$, where
$\mathcal O_\infty=k[\![t^{-1}]\!]$. The action
$\iota:\Gamma\to\Aut(T)$ has kernel $Z(\Gamma)$. Put
$\bar\Gamma=\iota(\Gamma)$. With $e_1,e_2$ the standard basis vectors,
the vertices $v_n=[\mathcal O_\infty t^n e_1+\mathcal O_\infty e_2]$
form the fundamental ray \cite{Serre} used as the standard lift in
Lemma~\ref{lem:normalization}. Their stabilizers are
\[
 G_0=S,\qquad
 G_n=\left\{\begin{pmatrix}a&f\\0&a^{-1}\end{pmatrix}:
       a\in k^\times,\ f\in A_{\le n}\right\}\quad(n\ge1).
\]
The edge $\varepsilon_n=[v_n,v_{n+1}]$ has stabilizer $B(k)$ if $n=0$
and $G_n$ if $n\ge1$. If $L$ preserves degrees, $\Phi_L$ preserves
each of these stabilizers. Thus the rule
$h_L(\gamma x)=\Phi_L(\gamma)x$ for $x=v_n,\varepsilon_n$ is
well defined, respects incidence and has inverse $h_{L^{-1}}$.
It gives a tree automorphism satisfying
$h_L\iota(\gamma)h_L^{-1}=\iota(\Phi_L(\gamma))$.
Thus $h_L$ lies in the normalizer $N_{\Aut(T)}(\bar\Gamma)$. For the
normalizer of $\GL(A)$ in $\GL(k(t))$, see \cite{MSQI}. We take subgroup
closures in the congruence topology of $\Gamma$.

\begin{figure}[htbp]
\centering
\begin{tikzpicture}[font=\small,vertex/.style={circle,fill,inner sep=1.7pt}]
\node[vertex,label=left:{$1$}] (v) at (0,0) {};
\foreach \i/\y in {1/.86,2/0,3/-.86}{
 \node[vertex] (a\i) at (2.1,\y) {};
 \node[vertex] (b\i) at (4.4,\y) {};
 \node[vertex] (c\i) at (6.7,\y) {};
 \draw (v)--(a\i)--(b\i)--(c\i);
 \draw[densely dotted] (c\i)--(7.55,\y);
}
\node[above=3pt] at (a1) {$(k,+)$};
\node[above=3pt] at (b1) {$(k,+)^2$};
\node[above=3pt] at (c1) {$(k,+)^3$};
\node[anchor=west] at (7.9,0) {$q+1$ rays};
\end{tikzpicture}
\caption{Three of the $q+1$ rays of $H_\ell\backslash T$. The top labels are vertex
groups. Along every ray the edge group has index one at the shallower
vertex and $q$ at the next vertex.}
\label{fig:quotient}
\end{figure}
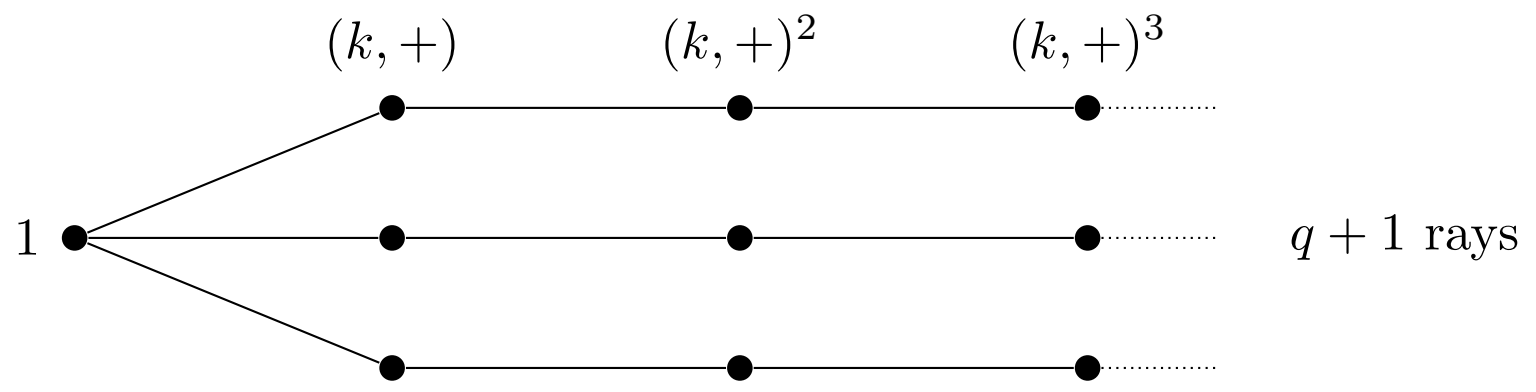

\begin{proposition}\label{prop:family}
The family $\{H_\ell\}$ is the full $\Aut(\Gamma)$-orbit of $\Gamma(t)$.
Every $H_\ell$ has index $|S|$ and acts faithfully on $T$. Its image
$\iota(H_\ell)$ is conjugate to $\iota(\Gamma(t))$ in $N_{\Aut(T)}(\bar\Gamma)$. The common
quotient graph of groups has one vertex with trivial group at depth zero
and $q+1$ rays, with vertex group $(k,+)^n$ at depth $n\ge1$.
Each edge group equals the group at its shallower endpoint and has index
$q$ at its deeper endpoint. The congruence members are precisely
$H_{\ev_c}=\Gamma(t-c)$ for $c\in k$. For $q\ge4$ every other $H_\ell$
is congruence-dense.
\end{proposition}
\begin{proof}
Define $L(1)=1$ and $L(t^n)=t^n-\ell(t^n)$ for $n\ge1$. The inverse
adds $\ell(t^n)$, so $L$ preserves degrees and
$\ell=\ev_0\circ L^{-1}$. On the amalgam factors this gives
$r_\ell=r_{\ev_0}\circ\Phi_{L^{-1}}$. Hence $H_\ell=\Phi_L(\Gamma(t))$
and $h_L$ gives the conjugacy. The retraction fixes $S$, so $H_\ell$
has index $|S|$ and intersects $Z(\Gamma)$ trivially.

Conversely the family is invariant under every automorphism in
Lemma~\ref{lem:normalization}. The amalgam gives
$\Phi_M(H_\ell)=H_{\ell\circ M^{-1}}$ for any $M\in\mathscr R_q$.
An entrywise ring automorphism $\sigma_*$ sends $H_\ell$ to
$H_{\tau\circ\ell\circ\sigma^{-1}}$ with $\tau=\sigma|_k$.
For $g\in\GL(A)$ write $g=\gamma\diag(\det g,1)$ with $\gamma\in\Gamma$.
Inner conjugations fix $H_\ell$ by normality. Constant diagonal conjugation commutes with $r_\ell$ because
$\ell((\det g)f)=(\det g)\ell(f)$. This proves the orbit assertion.

The vertex and edge groups are the stabilizers of chosen lifts. For
$\Gamma(t)$ the vertices over $v_n$ are indexed by $S/\rho_t(G_n)$,
a singleton at $n=0$ and $S/B(k)$ at $n\ge1$. Edges over each
$\varepsilon_n$ have the same labels $S/B(k)$. At positive depth the
endpoint maps preserve the labels, giving $q+1$ rays. The stabilizer is
trivial at $v_0$ and equals $\Gamma(t)\cap G_n=u(tA_{\le n-1})\simeq(k,+)^n$
at $v_n$ for $n\ge1$. Each edge stabilizer is the group at its shallower
endpoint (Figure~\ref{fig:quotient}) \cite{Bass93,Serre}.

If $H_\ell$ is congruence, Lemma~\ref{lem:kernel}(ii) gives
$H_\ell=\Gamma(t-c)$ for some $c\in k$. Its quasi-level is $(t-c)A$,
so $\ell(f)=f(c)$ since $\ell$ fixes constants. Conversely,
$H_{\ev_c}=\Gamma(t-c)$ is congruence. For $q\ge4$ apply
Corollary~\ref{cor:density} with $d=1$ to the other members.
\end{proof}

For $q=2,3$, let $\eta:S\to(k,+)$ be the abelianization map with
$\eta(u(1))=1$, and put $\chi_\ell=\eta\circ r_\ell$. A \emph{character}
is a homomorphism $\Gamma\to(k,+)$, called \emph{congruence} if it
factors through some $\rho_a$. Write $\mathcal V_q=\{\ev_c:c\in k\}$
and $D_c(f)=f'(c)$, using the formal derivative. The forms $\ell$
for which $\chi_\ell$ is congruence are exactly
\begin{align*}
 \mathcal E_2&=\{a_0\ev_0+a_1\ev_1+b_0D_0+b_1D_1:
                     a_i,b_i\in\F_2,\ a_0+a_1=1\},\\
 \mathcal E_3&=\bigl\{\textstyle\sum_{c\in\F_3}a_c\ev_c:
                     a_c\in\F_3,\ \textstyle\sum_c a_c=1\bigr\}.
\end{align*}

\begin{proposition}[Congruence closures for $q=2,3$]\label{prop:small}
For each $q\in\{2,3\}$ the congruence closure is
\[
 \cl(H_\ell)=
 \begin{cases}
   H_\ell=\Gamma(t-c), &\ell=\ev_c,\\
   \ker\chi_\ell, &\ell\in\mathcal E_q\setminus\mathcal V_q,\\
   \Gamma, &\ell\notin\mathcal E_q .
 \end{cases}
\]
Exactly six noncongruence members have proper closure. Each of these
closures has index $q$ in $\Gamma$. For each such $H_\ell$, the quotient
$\cl(H_\ell)/H_\ell$ is cyclic of order three for $q=2$ and quaternion
of order eight for $q=3$.
\end{proposition}
\begin{proof}
For $q=2$ the proper normal subgroups of $S\simeq S_3$ are $1$ and
$A_3=\ker\eta$. For $q=3$ those of $S\simeq Q_8\rtimes C_3$ are
$1$, $Z(S)$ and $Q_8=\ker\eta$. Suppose that $\ell$ is not an evaluation.
If $r_\ell(\Gamma(a))$ were central for some $a\ne0$,
Lemma~\ref{lem:kernel}(ii) would give $H_\ell=\Gamma(t-c)$, forcing
$\ell=\ev_c$. Thus if $H_\ell\Gamma(a)$ is a proper subgroup of $\Gamma$,
it equals $r_\ell^{-1}(\ker\eta)=\ker\chi_\ell$. Conversely, if
$\chi_\ell$ factors through $\rho_a$, then $H_\ell\Gamma(a)\subseteq
\ker\chi_\ell$ and equality follows. Intersecting these products gives
$\ker\chi_\ell$ when $\chi_\ell$ is congruence and $\Gamma$ otherwise.

Since the conjugates of $u(A)$ generate $\Gamma$, a character is determined
by its restriction to $u(A)$, where $\chi_\ell(u(f))=\ell(f)$. Identify the
complex $q$th roots of unity with $(k,+)$. By Boylan and Skoruppa
\cite[Theorem~1]{BS}, for $k=\F_3$ the restrictions of congruence characters
are exactly the linear combinations of the three evaluations. For $k=\F_2$,
each degree-one prime $z=t-c$ contributes the evaluation at $c$ and an
independent character factoring through reduction modulo $z^2$. Put $x=f(c)$ and
$y=f'(c)$. Modulo $z^2$ we have $f(c+z)=x+zy$ and
$u(x+zy)=u(x)(I+zyE_{12})$, with $E_{12}$ the upper matrix unit.
By \cite[(2)]{BS}, the extra character takes the value $(-1)^y$ on
$u(x+zy)$, or $y$ in additive notation. Its value on $u(z)$ is nonzero, so its
exact level is $z^2$. This proves the description of $\mathcal E_q$.

For $q=3$ the three evaluations are independent by interpolation. For
$q=2$ the two values and two derivatives are independent by Chinese
remaindering modulo $t^2(t-1)^2$. Hence $|\mathcal E_3|=3^2$ and
$|\mathcal E_2|=2^3$. Removing the $q$ evaluations leaves six members
in each case. Since $r_\ell$ is onto, $\Gamma/\ker\chi_\ell\simeq(k,+)$
and $\ker\chi_\ell/H_\ell\simeq\ker\eta$.
\end{proof}

For instance, over $\F_2$ the form $\ell=\ev_0+D_0$ takes the values
$1,0$ on $t,t^2$, so it differs from both evaluations. The subgroup
$H_\ell$ has index six, and its index-two closure contains $\Gamma(t^2)$.
Over $\F_3$, the form $\ell=2\ev_0+2\ev_1$ takes the value $2$ on
both $t$ and $t^2$, so it differs from every evaluation. Here $H_\ell$
has index $24$ and its closure has index three.

The automorphic images of each $\Gamma(t-c)$, $c\in k$, form this same
family by Proposition~\ref{prop:family}. The values $\ell(t^n)$, $n\ge1$,
are arbitrary. Distinct choices give distinct subgroups because $\ker\ell$
and $\ell(1)=1$ determine $\ell$. Only finitely many members have proper
closure, so for every $q$ there are continuum many congruence-dense
members with conjugate actions on $T$.

For $H=H_\ell$, write $[v]=Hv$ for the orbit of a vertex $v$ and
$V(H\backslash T)$ for the vertex set of the quotient graph.
Let $H_x=\operatorname{Stab}_H(x)$ for each vertex or edge $x$ of $T$.
Give $[v]$ mass $\mu_H([v])=1/|H_v|$. The vertex-averaging operator
on $\ell^2(V(H\backslash T),\mu_H)$ is
\[
   (P_Hf)([v])=\frac1{q+1}\sum_{w\sim v}f([w]),
\]
where $w$ ranges over the neighbors of $v$ in $T$. A quotient edge with
lift $e=(v,w)$ has multiplicity $[H_v:H_e]$ at $[v]$, and
$\mu_H([v])[H_v:H_e]=1/|H_e|=\mu_H([w])[H_w:H_e]$ gives detailed
balance. The tree conjugacies preserve these masses and multiplicities,
so pullback unitarily intertwines $P_H$, weighted adjacency $(q+1)P_H$,
and the normalized Laplacian $I-P_H$.

\end{document}